\documentstyle[12pt]{article}
\newcommand{\CC}{{\Bbb C}}
\newcommand{\bd}{\bar{d}}
\newcommand{\GL}{{\rm GL}}

\newcommand{\M}{{\cal M}}
\newcommand{\Mnd}{{\cal M}(n,d)}
\newcommand{\n}{\bar{n}}

\newcommand{\cO}{{\cal O}}
\newcommand{\Pic}{{\rm Pic}}
\newcommand{\PP}{{\Bbb P}}

\newcommand{\qed}{\hfill $\Box$}
\newcommand{\rk}{{\rm \ rk\, }}

\newcommand{\SL}{{\rm SL}}
\newcommand{\SM}{{\cal SM}}
\newcommand{\SMnd}{{\cal SM}(n,d)}
\newcommand{\SMnL}{{\cal SM}(n,L)}
\newcommand{\ZZ}{{\Bbb Z}}
\newtheorem{thm}{Theorem}
\newtheorem{conj}[thm]{Conjecture}
\newtheorem{cor}[thm]{Corollary}

\newtheorem{lem}[thm]{Lemma}
\newtheorem{prop}[thm]{Proposition}

\begin{document}
\makeatletter

%

\catcode`\@=11

%
\ifcase\@ptsize \font\tenmsx=msxm10\or
                \font\tenmsx=msxm10 scaled\magstephalf\or
                \font\tenmsx=msxm10 scaled\magstep1\fi
\ifcase\@ptsize \font\sevenmsx=msxm7\or
                \font\sevenmsx=msxm7 scaled\magstephalf\or
                \font\sevenmsx=msxm7 scaled\magstep1\fi
\ifcase\@ptsize \font\fivemsx=msxm5\or
                \font\fivemsx=msxm5 scaled\magstephalf\or
                \font\fivemsx=msxm5 scaled\magstep1\fi
\ifcase\@ptsize \font\tenmsy=msym10\or
                \font\tenmsy=msym10 scaled\magstephalf\or
                \font\tenmsy=msym10 scaled\magstep1\fi
\ifcase\@ptsize \font\sevenmsy=msym7\or
                \font\sevenmsy=msym7 scaled\magstephalf\or
                \font\sevenmsy=msym7 scaled\magstep1\fi
\ifcase\@ptsize \font\fivemsy=msym5\or
                \font\fivemsy=msym5 scaled\magstephalf\or
                \font\fivemsy=msym5 scaled\magstep1\fi
\newfam\msxfam
\newfam\msyfam
\textfont\msxfam=\tenmsx  \scriptfont\msxfam=\sevenmsx
  \scriptscriptfont\msxfam=\fivemsx
\textfont\msyfam=\tenmsy  \scriptfont\msyfam=\sevenmsy
  \scriptscriptfont\msyfam=\fivemsy

\def\hexnumber@#1{\ifcase#1 0\or1\or2\or3\or4\or5\or6\or7\or8\or9\or
	A\or B\or C\or D\or E\or F\fi }

\expandafter\ifx\csname frak\endcsname\relax
\ifcase\@ptsize \font\teneuf=eufm10\or
                \font\teneuf=eufm10 scaled\magstephalf\or
                \font\teneuf=eufm10 scaled\magstep1\fi
\ifcase\@ptsize \font\seveneuf=eufm7\or
                \font\seveneuf=eufm7 scaled\magstephalf\or
                \font\seveneuf=eufm7 scaled\magstep1\fi
\ifcase\@ptsize \font\fiveeuf=eufm5\or
                \font\fiveeuf=eufm5 scaled\magstephalf\or
                \font\fiveeuf=eufm5 scaled\magstep1\fi
  \newfam\euffam
  \textfont\euffam=\teneuf
  \scriptfont\euffam=\seveneuf
  \scriptscriptfont\euffam=\fiveeuf
  \def\frak{\ifmmode\let\next\frak@\else
   \def\next{\errmessage{Use \string\frak\space only in math mode}}\fi\next}
  \def\goth{\ifmmode\let\next\frak@\else
   \def\next{\errmessage{Use \string\goth\space only in math mode}}\fi\next}
  \def\frak@#1{{\frak@@{#1}}}
  \def\frak@@#1{\fam\euffam#1}
\fi
\edef\msx@{\hexnumber@\msxfam}
\edef\msy@{\hexnumber@\msyfam}

\mathchardef\boxdot="2\msx@00
\mathchardef\boxplus="2\msx@01
\mathchardef\boxtimes="2\msx@02
\mathchardef\square="0\msx@03
\mathchardef\blacksquare="0\msx@04
\mathchardef\centerdot="2\msx@05
\mathchardef\lozenge="0\msx@06
\mathchardef\blacklozenge="0\msx@07
\mathchardef\circlearrowright="3\msx@08
\mathchardef\circlearrowleft="3\msx@09
\mathchardef\rightleftharpoons="3\msx@0A
\mathchardef\leftrightharpoons="3\msx@0B
\mathchardef\boxminus="2\msx@0C
\mathchardef\Vdash="3\msx@0D
\mathchardef\Vvdash="3\msx@0E
\mathchardef\vDash="3\msx@0F
\mathchardef\twoheadrightarrow="3\msx@10
\mathchardef\twoheadleftarrow="3\msx@11
\mathchardef\leftleftarrows="3\msx@12
\mathchardef\rightrightarrows="3\msx@13
\mathchardef\upuparrows="3\msx@14
\mathchardef\downdownarrows="3\msx@15
\mathchardef\upharpoonright="3\msx@16
\let\restriction=\upharpoonright
\mathchardef\downharpoonright="3\msx@17
\mathchardef\upharpoonleft="3\msx@18
\mathchardef\downharpoonleft="3\msx@19
\mathchardef\rightarrowtail="3\msx@1A
\mathchardef\leftarrowtail="3\msx@1B
\mathchardef\leftrightarrows="3\msx@1C
\mathchardef\rightleftarrows="3\msx@1D
\mathchardef\Lsh="3\msx@1E
\mathchardef\Rsh="3\msx@1F
\mathchardef\rightsquigarrow="3\msx@20
\mathchardef\leftrightsquigarrow="3\msx@21
\mathchardef\looparrowleft="3\msx@22
\mathchardef\looparrowright="3\msx@23
\mathchardef\circeq="3\msx@24
\mathchardef\succsim="3\msx@25
\mathchardef\gtrsim="3\msx@26
\mathchardef\gtrapprox="3\msx@27
\mathchardef\multimap="3\msx@28
\mathchardef\therefore="3\msx@29
\mathchardef\because="3\msx@2A
\mathchardef\doteqdot="3\msx@2B
\let\Doteq=\doteqdot
\mathchardef\triangleq="3\msx@2C
\mathchardef\precsim="3\msx@2D
\mathchardef\lesssim="3\msx@2E
\mathchardef\lessapprox="3\msx@2F
\mathchardef\eqslantless="3\msx@30
\mathchardef\eqslantgtr="3\msx@31
\mathchardef\curlyeqprec="3\msx@32
\mathchardef\curlyeqsucc="3\msx@33
\mathchardef\preccurlyeq="3\msx@34
\mathchardef\leqq="3\msx@35
\mathchardef\leqslant="3\msx@36
\mathchardef\lessgtr="3\msx@37
\mathchardef\backprime="0\msx@38
\mathchardef\risingdotseq="3\msx@3A
\mathchardef\fallingdotseq="3\msx@3B
\mathchardef\succcurlyeq="3\msx@3C
\mathchardef\geqq="3\msx@3D
\mathchardef\geqslant="3\msx@3E
\mathchardef\gtrless="3\msx@3F
\mathchardef\sqsubset="3\msx@40
\mathchardef\sqsupset="3\msx@41
\mathchardef\vartriangleright="3\msx@42
\mathchardef\vartriangleleft="3\msx@43
\mathchardef\trianglerighteq="3\msx@44
\mathchardef\trianglelefteq="3\msx@45
\mathchardef\bigstar="0\msx@46
\mathchardef\between="3\msx@47
\mathchardef\blacktriangledown="0\msx@48
\mathchardef\blacktriangleright="3\msx@49
\mathchardef\blacktriangleleft="3\msx@4A
\mathchardef\vartriangle="0\msx@4D
\mathchardef\blacktriangle="0\msx@4E
\mathchardef\triangledown="0\msx@4F
\mathchardef\eqcirc="3\msx@50
\mathchardef\lesseqgtr="3\msx@51
\mathchardef\gtreqless="3\msx@52
\mathchardef\lesseqqgtr="3\msx@53
\mathchardef\gtreqqless="3\msx@54
\mathchardef\Rrightarrow="3\msx@56
\mathchardef\Lleftarrow="3\msx@57
\mathchardef\veebar="2\msx@59
\mathchardef\barwedge="2\msx@5A
\mathchardef\doublebarwedge="2\msx@5B
\mathchardef\angle="0\msx@5C
\mathchardef\measuredangle="0\msx@5D
\mathchardef\sphericalangle="0\msx@5E
\mathchardef\varpropto="3\msx@5F
\mathchardef\smallsmile="3\msx@60
\mathchardef\smallfrown="3\msx@61
\mathchardef\Subset="3\msx@62
\mathchardef\Supset="3\msx@63
\mathchardef\Cup="2\msx@64
\let\doublecup=\Cup
\mathchardef\Cap="2\msx@65
\let\doublecap=\Cap
\mathchardef\curlywedge="2\msx@66
\mathchardef\curlyvee="2\msx@67
\mathchardef\leftthreetimes="2\msx@68
\mathchardef\rightthreetimes="2\msx@69
\mathchardef\subseteqq="3\msx@6A
\mathchardef\supseteqq="3\msx@6B
\mathchardef\bumpeq="3\msx@6C
\mathchardef\Bumpeq="3\msx@6D
\mathchardef\lll="3\msx@6E
\let\llless=\lll
\mathchardef\ggg="3\msx@6F
\let\gggtr=\ggg
\mathchardef\circledS="0\msx@73
\mathchardef\pitchfork="3\msx@74
\mathchardef\dotplus="2\msx@75
\mathchardef\backsim="3\msx@76
\mathchardef\backsimeq="3\msx@77
\mathchardef\complement="0\msx@7B
\mathchardef\intercal="2\msx@7C
\mathchardef\circledcirc="2\msx@7D
\mathchardef\circledast="2\msx@7E
\mathchardef\circleddash="2\msx@7F
\def\ulcorner{\delimiter"4\msx@70\msx@70 }
\def\urcorner{\delimiter"5\msx@71\msx@71 }
\def\llcorner{\delimiter"4\msx@78\msx@78 }
\def\lrcorner{\delimiter"5\msx@79\msx@79 }
\def\yen{\mathhexbox\msx@55 }
\def\checkmark{\mathhexbox\msx@58 }
\def\circledR{\mathhexbox\msx@72 }
\def\maltese{\mathhexbox\msx@7A }
\mathchardef\lvertneqq="3\msy@00
\mathchardef\gvertneqq="3\msy@01
\mathchardef\nleq="3\msy@02
\mathchardef\ngeq="3\msy@03
\mathchardef\nless="3\msy@04
\mathchardef\ngtr="3\msy@05
\mathchardef\nprec="3\msy@06
\mathchardef\nsucc="3\msy@07
\mathchardef\lneqq="3\msy@08
\mathchardef\gneqq="3\msy@09
\mathchardef\nleqslant="3\msy@0A
\mathchardef\ngeqslant="3\msy@0B
\mathchardef\lneq="3\msy@0C
\mathchardef\gneq="3\msy@0D
\mathchardef\npreceq="3\msy@0E
\mathchardef\nsucceq="3\msy@0F
\mathchardef\precnsim="3\msy@10
\mathchardef\succnsim="3\msy@11
\mathchardef\lnsim="3\msy@12
\mathchardef\gnsim="3\msy@13
\mathchardef\nleqq="3\msy@14
\mathchardef\ngeqq="3\msy@15
\mathchardef\precneqq="3\msy@16
\mathchardef\succneqq="3\msy@17
\mathchardef\precnapprox="3\msy@18
\mathchardef\succnapprox="3\msy@19
\mathchardef\lnapprox="3\msy@1A
\mathchardef\gnapprox="3\msy@1B
\mathchardef\nsim="3\msy@1C
\mathchardef\ncong="3\msy@1D
\def\napprox{\not\approx}
\mathchardef\varsubsetneq="3\msy@20
\mathchardef\varsupsetneq="3\msy@21
\mathchardef\nsubseteqq="3\msy@22
\mathchardef\nsupseteqq="3\msy@23
\mathchardef\subsetneqq="3\msy@24
\mathchardef\supsetneqq="3\msy@25
\mathchardef\varsubsetneqq="3\msy@26
\mathchardef\varsupsetneqq="3\msy@27
\mathchardef\subsetneq="3\msy@28
\mathchardef\supsetneq="3\msy@29
\mathchardef\nsubseteq="3\msy@2A
\mathchardef\nsupseteq="3\msy@2B
\mathchardef\nparallel="3\msy@2C
\mathchardef\nmid="3\msy@2D
\mathchardef\nshortmid="3\msy@2E
\mathchardef\nshortparallel="3\msy@2F
\mathchardef\nvdash="3\msy@30
\mathchardef\nVdash="3\msy@31
\mathchardef\nvDash="3\msy@32
\mathchardef\nVDash="3\msy@33
\mathchardef\ntrianglerighteq="3\msy@34
\mathchardef\ntrianglelefteq="3\msy@35
\mathchardef\ntriangleleft="3\msy@36
\mathchardef\ntriangleright="3\msy@37
\mathchardef\nleftarrow="3\msy@38
\mathchardef\nrightarrow="3\msy@39
\mathchardef\nLeftarrow="3\msy@3A
\mathchardef\nRightarrow="3\msy@3B
\mathchardef\nLeftrightarrow="3\msy@3C
\mathchardef\nleftrightarrow="3\msy@3D
\mathchardef\divideontimes="2\msy@3E
\mathchardef\varnothing="0\msy@3F
\mathchardef\nexists="0\msy@40
\mathchardef\mho="0\msy@66
\mathchardef\eth="0\msy@67
\mathchardef\eqsim="3\msy@68
\mathchardef\beth="0\msy@69
\mathchardef\gimel="0\msy@6A
\mathchardef\daleth="0\msy@6B
\mathchardef\lessdot="3\msy@6C
\mathchardef\gtrdot="3\msy@6D
\mathchardef\ltimes="2\msy@6E
\mathchardef\rtimes="2\msy@6F
\mathchardef\shortmid="3\msy@70
\mathchardef\shortparallel="3\msy@71
\mathchardef\smallsetminus="2\msy@72
\mathchardef\thicksim="3\msy@73
\mathchardef\thickapprox="3\msy@74
\mathchardef\approxeq="3\msy@75
\mathchardef\succapprox="3\msy@76
\mathchardef\precapprox="3\msy@77
\mathchardef\curvearrowleft="3\msy@78
\mathchardef\curvearrowright="3\msy@79
\mathchardef\digamma="0\msy@7A
\mathchardef\varkappa="0\msy@7B
\mathchardef\hslash="0\msy@7D
\mathchardef\hbar="0\msy@7E
\mathchardef\backepsilon="3\msy@7F

\def\Bbb{\ifmmode\let\next\Bbb@\else
 \def\next{\errmessage{Use \string\Bbb\space only in math mode}}\fi\next}
\def\Bbb@#1{{\Bbb@@{#1}}}
\def\Bbb@@#1{\fam\msyfam#1}

\catcode`\@=\active
\makeatother
\title{Theta Functions for $\SL(n)$ versus $\GL(n)$}
\author{Ron Donagi and Loring W. Tu}
\date{November 1, 1992}
\maketitle

\begin{enumerate}
\item[\S 1.] Theta bundles
\item[\S 2.] A Galois covering
\item[\S 3.] Pullbacks
\item[\S 4.] Proof of Theorem 1
\item[\S 5.] A conjectural duality
\end{enumerate}

Over a smooth complex projective curve $C$ of genus $g$
one may consider two types of
moduli spaces of vector bundles, $\M :=\Mnd$, the moduli space of semistable
bundles of rank $n$ and degree $d$ on $C$, and $\SM := \SMnL$, the moduli
space of those bundles whose determinant is isomorphic to a fixed
line bundle $L$ on $C$.
We call the former a {\it full moduli space} and the latter a
{\it fixed-determinant moduli space}.
Since the spaces $\SM (n,L)$ are all isomorphic as $L$ varies in $\Pic^d(C)$,
we also write $\SM (n,d)$ to denote any one of them.

On both moduli spaces there are well-defined {\it theta bundles},
as we recall in Section \ref{bundles}.
While the theta bundle $\theta$ on $\SM$ is uniquely defined, the theta
bundles $\theta _F$ on $\M$ depend on the choice of complementary vector
bundles $F$ of minimal rank over $C$.
For any positive integer $k$, sections of $\theta_F^k$ generalize the
classical theta functions of level $k$ on the Jacobian of a curve, and so we
call sections of $\theta^k$ over $\SM$ and $\theta_F^k$ over $\M$ {\it theta
functions of level $k$ for} $\SL (n)$ {\it and} $\GL (n)$ respectively.

Our goal is to study the relationship between these two spaces of theta
functions.
We prove a simple formula relating their dimensions, and then formulate a
conjectural duality between these two spaces.

\begin{thm}
	\label{formula}
If $h= \gcd (n,d)$ is the greatest common divisor of $n$ and $d$, and
$L \in \Pic ^d (C)$, then
$$\dim H^0(\SM (n,L), \theta^k )\cdot k^g =
\dim H^0(\M (n,d), \theta_F^k )\cdot h^g.$$
\end{thm}

Faltings \cite{faltings} has proven the Verlinde formula for semisimple
groups, which gives in particular the dimension of $H^0(\SM, \theta^k)$.
The dimension of $H^0(\M,\theta_F^k)$ is thus determined by Theorem
\ref{formula}.
When $k=1$ and $d=0$, \cite{beauville-narasimhan-ramanan} computes explicitly
the two spaces in Theorem \ref{formula}.
Their result is a forerunner of Theorem \ref{formula}.

Theorem \ref{formula} is consistent with and therefore lends credence to
another, so far conjectural,
relationship between these two types of theta functions.
To explain this, start with integers $\n, \bd, h, k$ such that $\n, h, k$ are
positive and $\gcd (\n,\bd) = 1$.
Let $F \in \M (\n,\bd)$ and write
$$\SM_1 = \SM (h \n, (\det F)^h)\quad { \rm and }\quad
\M_2 = \M (k\n, k(\n (g-1)-\bd)).$$
The tensor product map $\tau$ sends $\SM_1 \times \M_2$ to
$\M (hk\n ^2, hk\n^2 (g-1))$.

\begin{conj}
	\label{conjecture}
The tensor product map induces a natural duality between
$H^0(\SM _1, \theta^k)$ and \linebreak $H^0(\M_2, \theta_F^h)$.
\end{conj}

For further discussion of this duality, including supporting evidence, see
Section \ref{duality}.

\bigskip
\leftline{{\bf Notation and Conventions.}}
\medskip
\begin{enumerate}
\item[] $h^0(\ ) = \dim H^0(\ )$

\item[] $J_d = \Pic^d(C) = $\{isomorphism classes of line bundles of degree
$d$ on $C \}$

\item[] $J = J_0 = \Pic^0(C) $

\item[] $L_1 \boxtimes L_2 = \pi_1^*L_1 \otimes \pi_2^*L_2$ if $L_i$ is a
line bundle on $X_i$ and $\pi_i : X_1 \times X_2 \to X_i$ is the $i$-th
projection

\item[] $S^hC =$ the $h$th symmetric product of $C$

\item[] $T_n =$ the group of $n$-torsion bundles on $C$

\end{enumerate}

\section{Theta bundles}
	\label{bundles}

We recall here
the definitions of the theta bundles on a fixed-determinant moduli
space and on a full moduli space.
Our definitions are slightly different from but equivalent to those in
\cite{drezet-narasimhan}.

For $L \in \Pic ^d (C)$,  the Picard group of $\SM := \SM (n,L)$ is $\ZZ$ and
the theta bundle $\theta$ on $\SM$ is the positive generator of $\Pic (\SM)$.

When $n$ and $d$ are such that $\chi (E) = 0$ for $E \in \M (n,d)$,
i. e., when $d=(g-1)n$, there is a natural divisor $\Theta \subset \M
(n,n(g-1)) $:
$$\Theta = \mbox{closure of } \{ E \mbox{ stable in } \M(n,n(g-1))\ | \
h^0 (E) \ne 0 \}.$$
The theta bundle $\theta$ over $\M (n, n(g-1))$
is the line bundle corresponding to this divisor.

We say that a semistable
bundle $F$ is {\it complementary} to another bundle $E$ if
$\chi (E\otimes F) = 0$.
We also say that $F$ is {\it complementary} to $\M (n,d)$ if $\chi (E \otimes
F) = 0$ for any $E \in \M (n,d)$.
It follows easily from the Riemann-Roch theorem that if $E \in \M (n,d)$,
$h=\gcd (n,d), n= h \n$, and $d= h\bd$, then $F$ has rank $n_F$ and
degree $d_F$, where
$$n_F = k\n  \quad \mbox{ and }\quad d_F = k(\n (g-1) - \bd)$$
for some positive integer $k$.

If $F$ is complementary to $\M(n,d)$, let
$$\tau _F : \M (n,d) \to \M (nn_F, nn_F (g-1))$$
be the map
$$E \mapsto E \otimes F.$$
Pulling back the theta bundle $\theta$ from $\M (nn_F, nn_F (g-1))$ via
$\tau_F$ gives a line bundle $\theta_F := \tau_F^* \theta$ over $\M(n,d)$.
(This bundle may or may not correspond to a divisor in $\M(n,d)$.)
Let $\det: \M(n,d) \to J_d(C)$ be the determinant map.
When $\rk F$ is the minimal possible: $\rk F = \n =n/h$,
then $\theta_F$ is called a {\it theta bundle over} $\M (n,d)$;
otherwise, it is a multiple of a theta bundle.
Indeed, we extract from \cite{drezet-narasimhan} the formula:

\begin{prop}
	\label{thetaF}
Let $F$ and $F_0$ be two bundles complementary to $\M (n,d)$.
If $\rk F = a \rk F_0$, then
$$\theta_F \simeq \theta_{F_0}^{\otimes a} \otimes \det{}^*(\det F \otimes
(\det F_0)^{-a}),$$
where we employ the usual identification of $\Pic^0(C)$ with $\Pic^0 (J_0)$.
\end{prop}

In particular, $\theta_F$ depends only on $\rk F$ and $\det F$.

If $\theta_{F}$ is a theta bundle on $\M (n,d)$, then for any $L \in
\Pic^d(C)$, $\theta_{F}$ restricts to the theta bundle on $\SM (n,L)$.

\section{A Galois covering}
	\label{covering}

Let $\tau : Y \to X$ be a covering of varieties, by which we mean a finite
\'etale morphism.
A {\it deck transformation} of the covering is an automorphism $\phi : Y \to
Y$ that commutes with $\tau$.
The covering is said to be {\it Galois} if the group of deck transformations
acts transitively (hence simply transitively) on a general fiber
of the covering.

Denote by $J= \Pic ^0 (C)$ the group of isomorphism classes of line bundles
of degree 0 on the curve $C$, and $G = T_n$ the subgroup of torsion points of
order $n$.
Fix $L \in \Pic^d(C)$ and let $\SM = \SM (n,L)$, $J=J_0(C)$, and $\M =
\M(n,d)$.
Recall that the tensor product map
\begin{eqnarray*}
\tau : \SM \times J &\to& \M\\
(E, M) &\mapsto& E \otimes M
\end{eqnarray*}
gives an $n^{2g}$-sheeted \'etale map (\cite{teixidor-tu}, Prop. 8).
The group $G=T_n$ acts on $\SM \times J$ by
$$N.(E, M) = (E\otimes N^{-1}, N\otimes M).$$
It is easy to see that $G$ is the group of deck transformations of the
covering $\tau$ and that it acts transitively on every fiber.
Therefore, $\tau: \SM \times J \to \M$ is a Galois covering.

\begin{prop}
	\label{directimage}
If $\tau: Y \to X$ is a Galois covering with finite abelian Galois
group $G$, then $\tau_*\cO_Y$ is a vector bundle on $X$ which decomposes into
a direct sum of line bundles indexed by the characters of $G$:
$$\tau_*\cO_Y = \sum_{\lambda \in \hat G} L_{\lambda},$$
where $\hat G$ is the character group of $G$.
\end{prop}

\medskip
\noindent
{\sc Proof}.  Write $\cO= \cO_Y$.  The fiber of $\tau_*\cO$ at a point $x\in X$
is naturally a complex vector space with basis $\tau^{-1}(x)$.
Hence, $\tau_*\cO$ is a vector bundle over $X$.
The action of $G$ on $\tau^{-1}(x)$ induces a representation of $G$ on
$(\tau_*\cO)(x)$ equivalent to the regular representation.
Because $G$ is a finite abelian group, this representation of $G$ decomposes
into a direct sum of one-dimensional representations indexed by the
characters of $G$:
$$(\tau_*\cO)(x) = \sum_{\lambda \in \hat G} L_{\lambda}(x).$$
Thus, for every $\lambda \in \hat G$, we obtain a line bundle $L_\lambda$ on
$X$ such such $\tau_*\cO = \sum_{\lambda} L_\lambda$.  \quad\quad$\Box$

\section{Pullbacks}
	\label{pullbacks}

We consider the tensor product map
\begin{eqnarray*}
\tau: \SM (n_1, L_1) \times \M (n_2, d_2) &\to& \M(n_1n_2, n_1d_2+n_2d_1)\\
(E_1, E_2) &\mapsto& E_1 \otimes E_2,
\end{eqnarray*}
where $d_1 = \deg L_1$.
For simplicity, in this section we write $\SM_1 = \SM (n_1, L_1)$,
$\M_2 = \M (n_2, d_2)$, and $\M_{12} = \M(n_1n_2, n_1d_2 + n_2d_1)$.

\begin{prop}
	\label{pullbackform}
Let $F=F_{12}$ be a bundle on $C$ complementary to $\M_{12}$.
Then
$$\tau^*\theta_F \simeq \theta^c \boxtimes \theta_{E_1 \otimes F}$$
for any $E_1 \in \SM (n_1, L_1)$, where
$$c:={ {n_2 \rk F}\over {\rk F_1}} = {{n_2\rk F} \over {n_1/\gcd
(n_1,d_1)}}$$
and $F_1$ is a minimal complementary bundle to $E_1$.
\end{prop}

\medskip
\noindent
{\bf Proof.}  For $E_2 \in \M(n_2,d_2)$, let
$$\tau_{E_2} : \SM_1 \to \M_{12}$$
be tensoring with $E_2$.
Then
$$(\tau^*\theta_F)|_{\SM \times \{E_2\}} = \tau_{E_2}^*\theta_F =
\tau_{E_2}^*\tau_F^*\theta = \tau_{E_2 \otimes F}^* \theta = \theta^c,$$
where by Proposition \ref{thetaF}
\begin{eqnarray*}
c&=& \rk (E_2 \otimes F)/ \rk F_1 \\
&=&  {{n_2 \rk F}\over {n_1/\gcd(n_1,d_1)}}.
\end{eqnarray*}
Similarly,
\begin{eqnarray*}
(\tau^*\theta_F)|_{\{E_1\}\times\M_2} &=& \tau_{E_1}^*\theta_F =
\tau_{E_1}^*\tau_F^*\theta \\
&=& \tau_{E_1\otimes F}^*\theta = \theta_{E_1\otimes F}.
\end{eqnarray*}
Note that the bundle $\theta_{E_1\otimes F}$ depends only on
$\rk (E_1\otimes F) = n_1\rk F$ and
$\det (E_1 \otimes F) = L_1^{\rk F} \otimes (\det F)^{n_1}$.
Hence, both $(\tau^*\theta_F)|_{\SM_1 \times\{E_2\}}$ and
$(\tau^*\theta_F)|_{\{E_1\} \times \M_2}$ are independent of $E_1$ and $E_2$.
By the seesaw theorem,
$$\tau^*\theta_F \simeq \theta^c \boxtimes \theta_{E_1\otimes F}.$$
\rightline{$\Box$}

\begin{cor}
	\label{jacobian}
Let $L \in \Pic^d(C)$ and
$$\tau: \SM (n,L) \times J_0 \to \M(n,d)$$
be the tensor product map.
Suppose $F$ is a minimal complementary bundle to $\M (n,d)$.
Choose $N \in \Pic^{g-1}(C)$ to be a line bundle
such that $N^n = L \otimes (\det F)^h$,
where $h=\gcd (n,d)$.
Then
$$\tau^* \theta_F = \theta \boxtimes \theta_N^{n^2/h}.$$
\end{cor}

\medskip
\noindent
{\bf Proof.}  Apply the Proposition with $\rk F = n/h$ and
$n_1 = n, d_1=d, n_2=1, d_2=0$.
Then $c=1$.
By Proposition \ref{thetaF},
\begin{eqnarray*}
\theta_{E_1\otimes F} &=& \theta_N^{n^2/h} \otimes \det{}^*
(\det (E_1 \otimes F)\otimes N^{-n^2/h}) \\
&=& \theta_N^{n^2/h}.
\end{eqnarray*}
\rightline{$\Box$}

\section{Proof of Theorem 1}
	\label{proof}

We apply the Leray spectral sequence to compute the cohomology of
$\tau^*\theta_F^k$ on the total space of the covering
$\tau : \SM \times J \to \M$ of Section \ref{covering}.
Recall that $\SM = \SM (n,d)$, $J=J_0$, and $\M = \M(n,d)$.
Because the fibers of $\tau$ are 0-dimensional, the spectral sequence
degenerates at the $E_2$-term and we have
\begin{equation}
	\label{a}
H^0(\SM \times J , \tau^*\theta_F^k ) = H^0(\M, \tau_*\tau^*\theta_F^k).
\end{equation}
By Cor. \ref{jacobian} and the K\"unneth formula, the left-hand side of
(\ref{a}) is
\begin{eqnarray*}
H^0(\SM \times J, \tau^*\theta_F^k))&=& H^0(\SM \times J, \theta^k \boxtimes
\theta_N^{kn^2/h}))\\
&=& H^0(\SM, \theta^k)\otimes H^0(J, \theta_N^{kn^2/h}).
\end{eqnarray*}
By the Riemann-Roch theorem for an abelian variety,
$$h^0(J, \theta_N^{kn^2/h}) = (kn^2/h)^g.$$
So the left-hand side of (\ref{a}) has dimension
\begin{equation}
	\label{left}
h^0(\SM, \theta^k) \cdot (kn^2/h)^g.
\end{equation}

Next we look at the right-hand side of (\ref{a}).
By the projection formula and Prop. \ref{directimage},
\begin{eqnarray*}
\tau_*\tau^*\theta_F^k &=& \theta_F^k \otimes \tau_*\cO \\
&=& \theta_F^k \otimes \sum_{\lambda \in {\hat G}} L_{\lambda}\\
&=& \sum_{\lambda \in {\hat G}} \theta_F^k \otimes L_{\lambda}.
\end{eqnarray*}

Our goal now is to show that for any character $\lambda \in \hat G$,
\begin{equation}
	\label{b}
H^0(\M, \theta_F^k \otimes L_\lambda) \simeq H^0 (\M, \theta_F^k).
\end{equation}
This will follow from two lemmas.

\begin{lem}
	\label{L}
The line bundle $L_{\lambda}$ on $\M$ is the pullback under $\det : \M \to
J_d$ of some line bundle $N_{\lambda}$ of degree 0 on $J_d := \Pic ^d(C)$.
\end{lem}

\begin{lem}
	\label{independence}
For $F$ a vector bundle as above, $k$ a positive integer, and $M$ a line
bundle of degree 0 over $C$,
$$H^0 (\M, \theta_{F\otimes M}^k) \simeq H^0(\M, \theta_F^k ).$$
\end{lem}

Assuming these two lemmas, let's prove (\ref{b}).
By Proposition \ref{thetaF},
$$\theta_{F\otimes M} = \theta_F \otimes \det {}^* M^{n_F};$$
hence,
$$\theta_{F\otimes M}^k = \theta_F^k \otimes \det {}^* M^{n_Fk}.$$
If $L_{\lambda}= \det{}^*N_\lambda$, and we choose a root
$M= N_{\lambda}^{1/(n_Fk)}$,
then
$$\theta_F^k \otimes L_{\lambda} = \theta_F^k \otimes \det{}^*N_{\lambda} =
\theta_{F\otimes M}^k.$$
Equation (\ref{b}) then follows from Lemma \ref{independence}.

\medskip
\noindent
{\sc Proof of Lemma} \ref{L}.
Define $\alpha : \SM \times J \to J$
to be the projection onto the second factor,
$\beta : \M \to J$
to be the composite of $\det : \M \to J_d$ followed by multiplication by
$L^{-1}: J_d \to J$, and
$\rho: J \to J$
to be the $n$-th tensor power map.
Then there is a commutative diagram
$$\begin{array}{rcccl}
\;     & \SM \times J & \stackrel{\tau}{\to} & \M        &\;    \\
       &              &                      &           &      \\
\alpha & \downarrow   & \;                   &\downarrow &\beta \\
       &              &                      &           &      \\
\;     & J            & \stackrel{\rho}{\to} &  J.       &\;
\end{array}$$

Furthermore, in the map $\alpha: \SM \times J \to J$ we let $G= T_n$ act on
$J$ by
$$N.M = N \otimes M, \quad\quad M \in J,$$
and in the map $\beta: \M \to J$ we let $G$ act trivially on both $\M$ and
$J$.
Then all the maps in the commutative diagram above are $G$-morphisms.

By the push-pull formula (\cite{hartshorne}, Ch. III, Prop. 9.3, p. 255),
$$ \tau_*\alpha^*\cO_J = \beta^*\rho_*\cO_J .$$
By Proposition \ref{directimage}, $\rho_*\cO_J$ is a direct sum of line
bundles $V_\lambda$ on $J$, where $\lambda \in \hat G$.
In fact, these $V_\lambda$ are precisely the $n-$torsion bundles in $J$;
in particular, their degrees are zero.
If $\tau_{L^{-1}} : J_d \to J$ is multiplication by the line bundle $L^{-1}$,
we set $N_\lambda := \tau_{L^{-1}}^* V_\lambda$.
Then
\begin{eqnarray*}
\tau_*\cO_{\SM \times J}
&=& \beta^*\sum_{\lambda \in \hat G} V_{\lambda} \\
&=& \det{}^*\tau_{L^{-1}}^* \sum V_\lambda \\
&=& \sum \det{}^*N_{\lambda}.
\end{eqnarray*}

By Prop. \ref{directimage}, $\tau_*\cO_{\SM \times J} = \sum L_\lambda$.
Since both $L_\lambda$ and $\det{}^* N_\lambda$ are eigenbundles of
$\tau_*\cO_{\SM \times J}$ corresponding to the character $\lambda \in \hat
G$,
$$L_\lambda = \det{}^* N_\lambda.$$
\qed

\medskip
\noindent
{\sc Proof of Lemma} \ref{independence}. Tensoring with $M\in J_0(C)$ gives an
automorphism
\begin{eqnarray*}
\tau_M : \M &\to& \M \\
E &\mapsto& E\otimes M,
\end{eqnarray*}
under which
$$\theta_{F\otimes M} = \tau_M^*\theta_F .$$
Hence,
$$\theta_{F\otimes M}^k = \tau_M^*(\theta_F^k)$$
and the lemma follows. \quad\quad\quad $\Box$

\medskip

Returning now to Eq. (\ref{a}), its right-hand side is
\begin{eqnarray*}
H^0(\M, \tau_*\tau^*\theta_F^k) &=& \sum _{\lambda \in \hat G} H^0(\M,
\theta_F^k \otimes L_\lambda ) \\
&\simeq& \sum _{\lambda \in \hat G} H^0 (\M, \theta_F^k),\quad\quad\mbox{ (by
(\ref{b}))}
\end{eqnarray*}
which has dimension
$$h^0(\M ,\theta_F^k) \cdot n^{2g}.$$
By (\ref{left}) the left-hand side of Eq. (\ref{a})
has dimension
$$h^0(\SM, \theta^k)\cdot (kn^2/h)^{g}.$$
Equating these two expressions gives
$$h^0(\M, \theta_F^k) = h^0(\SM,\theta^k) \cdot
({k\over h})^g.$$
This completes the proof of Theorem \ref{formula}.

\section{A conjectural duality}
	\label{duality}

As in the Introduction we start with integers $\n$, $\bd$, $h$, $k$ such
that $\n$, $h$, $k$ are positive and $\gcd (\n, \bd)=1$.
Take
$$n_1=h\n ,\ d_1=h\bd,\ n_2=k\n,\ d_2=k(\n(g-1)-\bd), \ {\rm and }\
L_1\in \Pic^{d_1}(C).$$
The tensor product induces a map
$$\tau: \SM (n_1, L_1) \times \M(n_2, d_2) \to \M (n_1n_2, n_1n_2(g-1)).$$
As before, write $\SM_1 = \SM (n_1, L_1)$, $\M_2= \M(n_2,d_2)$,
and $\M_{12} =\M(n_1n_2, n_1n_2(g-1))$.
Let $F_2=F$ and $F_{12}=\cO$ be complementary to $\M_2$ and $\M_{12}$
respectively.

By the pullback formula (Proposition \ref{pullbackform})
$$\tau^*\theta_{\cO} = \theta^{{n_2}/\n} \boxtimes \theta_{E_1}.$$
But by Proposition \ref{thetaF},
$$\theta_{E_1} = \theta_F^h \otimes \det{}^*(L\otimes (\det F)^{-h}).$$
If $L= (\det F)^h$, then $\theta_{E_1} = \theta_F^h$ and
$$\tau^*\theta_{\cO} = \theta ^k \boxtimes \theta_F^h.$$
By the K\"unneth formula,
$$H^0(\SM_1 \times \M_2, \tau^*\theta_{\cO})= H^0(\SM_1, \theta^k)\otimes
H^0(\M_2, \theta_F^h).$$
In \cite{beauville-narasimhan-ramanan} it is shown that
up to a constant, $\theta_{\cO}$ has a unique section $s$ over $\M_{12}$.
Then $\tau^*s$ is a section of $H^0(\SM_1 \times \M_2 , \tau^*\theta_{\cO})$
and therefore induces a natural map
\begin{equation}
	\label{dual}
H^0(\SM_1, \theta^k)^{\vee} \to H^0(\M_2, \theta_F^h).
\end{equation}
We conjecture that this natural map is an isomorphism.


Among the evidence for the duality (\ref{dual}), we cite the following.
\begin{enumerate}
\item[i)] (Rank 1 bundles) The results of
\cite{beauville-narasimhan-ramanan} that
$$H^0(\SM (n,\cO), \theta )^{\vee} \simeq H^0(\M (1,g-1), \theta_{\cO}^n)
\quad {\rm and } \quad H^0(\M (n,n(g-1)), \theta_{\cO} ) = \CC,$$
are special cases of (\ref{dual}), for $(n_2,d_2)=(1,g-1)$ and $(n_1,d_1)
=(1,0)$ respectively.

\item[ii)] (Consistency with Theorem \ref{formula}) Given a triple of
integers $(n_1, d_1, k)$, we define $h, \n, \bd$ by
$$h=\gcd (n_1,d_1), n_1=h\n, d_1=h\bd$$
and let $n_2, d_2$ be as before:
$$n_2=k\n, d_2=k(\n(g-1)-\bd).$$
Assuming $n_1$ and $k$ to be positive, it is easy to check that the function
$$(n_1,d_1,k) \mapsto (n_2, d_2, h)$$
is an involution.
Write $v(n,d,k)=h^0(\Mnd, \theta_F^k)$ and
$s(n,d,k)=h^0(\SMnd, \theta^k)$.
Then Theorem \ref{formula} assumes the form
\begin{equation}
	\label{one}
v(n,d,k) \cdot h^g = s(n,d,k)\cdot k^g.
\end{equation}
The duality (\ref{dual}) implies that
there is an equality of dimensions
\begin{equation}
	\label{two}
s(n_1,d_1, k) = v(n_2, d_2, h).
\end{equation}
Because $(n_1,d_1,k) \mapsto (n_2,d_2,h)$ is an involution,
it follows that
\begin{equation}
	\label{three}
s(n_2,d_2, h) = v(n_1, d_1, k).
\end{equation}
Putting (\ref{one}), (\ref{two}), and (\ref{three}) together, we get
$$v(n_2,d_2,h)k^g = s(n_2,d_2, h)h^g,$$
which is Theorem \ref{formula} again.

\item[iii)] (Elliptic curves) We keep the notation above, specialized to
the case of a curve $C$ of genus $g=1$:
$$n_1 = h\n, d_1 = h\bd, n_2= k\n, d_2= -k\bd.$$
Set $C' := \Pic^{\bd} (C) $.
The map sending a line bundle to its dual
gives an isomorphism $C'\simeq \Pic^{-\bd}(C)$.
If $L \in \Pic ^{\bd}(C)$, viewed as a line bundle on $C$, we let $\ell$
be the corresponding point in $C'$, and $\cO _{C'}(\ell)$ the associated line
bundle of degree 1 on the curve $C'$.
There is a natural map
$$\gamma : \Pic^{h\bd} (C) \to \Pic ^h (C')$$
which sends $L:= L_1\otimes \cdots \otimes L_h \in \Pic^{h\bd}(C)$ to
$L':= \cO_{C'} (\ell _1 + \cdots +\ell_h )$,
where $L_i \in \Pic^{\bd}(C)$ corresponds to the point $\ell _i \in C'$.

\end{enumerate}

 From \cite{atiyah} and \cite{tu} we see that there are natural identifications
$$\M(h\n, h\bd) \simeq S^h\M (\n,\bd) \simeq S^h \Pic^{\bd}(C) = S^hC'$$
and
$$\M(k\n, -k\bd) \simeq S^k\M (\n,-\bd) \simeq S^k \Pic^{-\bd}(C) \simeq
S^kC'.$$
Furthermore, there is a commutative diagram
$$\begin{array}{rcccl}
\;     & \M (h\n,h\bd) & \stackrel{\sim}{\to} & S^hC'        &\;    \\
       &               &                   &                 &      \\
\det & \downarrow   & \;                   &\downarrow &\alpha \\
       &               &                   &           &       \\
\;     & \Pic^{h\bd}(C)& \stackrel{\gamma}{\to} & \Pic^h(C').       &\;
\end{array}$$
Since the fiber of the Abel-Jacobi map $\alpha: S^hC' \to \Pic^h(C')$ above
$L'$ is the projective space $\PP H^0(C',L')$, it follows that there is a
natural identification
$$\SM (h\n, L) \simeq \PP H^0 (C', L').$$
Since the theta bundle is the positive generator of $\SM (h\n, L)$, it is
the hyperplane bundle.
For $F\in \M (\n, \bd)$, let $q\in C'$ be the point corresponding to the line
bundle $Q:= \det F \in \Pic^{\bd}(C)$.
Then
\begin{eqnarray*}
H^0(\SM(h\n, (\det F)^h), \theta^k) &\simeq& H^0 (\PP H^0(C', \cO_{C'}(hq)),
\cO(k))\\
&=& S^kH^0(C', \cO_{C'} (hq))^{\vee}.
\end{eqnarray*}

Recall that each point $q \in C'$ determines a divisor $X_q$ on the symmetric
product $S^kC'$:
$$X_q := \{ q+D \ | \ D\in S^{k-1}C'\}.$$
The proof of Theorem 6 in \cite{tu} actually shows that if $F \in
\M (\n,-\bd)$,
then under the identification $\M (k\n, -k\bd) \simeq S^kC'$,
the theta bundle $\theta_F$ corresponds to the bundle
associated to the divisor
$X_q$ on $S^kC'$, where $q$ is the point corresponding to $\det F \in
\Pic^{\bd}$.
Therefore, by the calculation of the cohomology of a symmetric product
in \cite{tu}
\begin{eqnarray*}
H^0 (\M (k\n, -k\bd), \theta_F^h) &=& H^0 (S^kC', \cO (h X_q))\\
&=& S^k H^0 (C', \cO (hq)).
\end{eqnarray*}
So the two spaces $H^0(\SM(h\n, (\det F)^h), \theta^k)$
and $H^0 (\M (k\n, -k\bd), \theta_F^h)$ are naturally dual to each other.

\begin{enumerate}

\item[iv)] (Degree 0 bundles) Consider the moduli space $\SM(n,0)$ of rank
$n$ and degree $0$ bundles.
In this case,
$$n_1=n,\ d_1=0,\ h=\gcd(n,0)=n,\ n_2=k,\ d_2=k(g-1).$$
So the conjectural duality is
$$H^0(\SM(n,\cO), \theta^k)^{\vee} \simeq H^0(\M(k,k(g-1)),\theta_{\cO}^n).$$
Because $\M (k,k(g-1))$ is isomorphic to $\M (k,0)$ (though noncanonically),
it follows that in the notation of ii)
$$s(n,0,k) = v(k,0, n).$$
According to R. Bott and A. Szenes, this equality follows from Verlinde's
formula.

\end{enumerate}

\end{document}